\documentclass[twoside,11pt]{article}
\usepackage{jnm}
\usepackage[T2A]{fontenc}
\usepackage{slashbox}
\usepackage{algorithmic} \usepackage[ruled]{algorithm}

\usepackage{graphicx}

\def\cond{\operatorname{cond}}
\def\eps{\varepsilon}
\def\diag{\operatorname{diag}}

\title{TT-GMRES: on solution to a linear system in the structured tensor format}

\author{Sergey V. Dolgov
\thanks{Institute of Numerical Mathematics, Russian Academy of Sciences, Gubkina Street, 8, Moscow, Russia}
\support{This work was supported in part by
         RFBR grants  09-01-12058, 10-01-00757, 11-01-00549, RFBR/DFG grant 09-01-91332,
         Russian Federation Gov. contracts No. $\Pi 1178,$ $\Pi 1112$ and $\Pi 940,$ $14.740.11.0345$
         and Promotionstipendium of Max-Planck Gesellschaft.
         Part of this work was done during the stay of the author in Max-Plank Institute for Mathematics in Sciences, Leipzig, Germany.}
}


\abstract{
A adapted tensor-structured GMRES method for the TT format is proposed and investigated. The Tensor Train (TT) approximation is a robust approach to high-dimensional problems. One class of problems is solution of a linear system. In this work we study the convergence of the GMRES method in the presence of tensor approximations and provide relaxation techniques to improve its performance. Several numerical examples are presented. The method is also compared with a projection TT linear solver based on the ALS and DMRG methods. On a particular sPDE (high-dimensional parametric) problem, these methods manifest comparable performance,
with a good preconditioner the TT-GMRES overcomes the ALS solver.
}

\begin{document}

\maketitle

\section{Introduction}
Solving linear systems arising from problems with many dimensions is a computationally demanding task.
Such problems are posed, for example, in quantum chemistry \cite{meyer-mctdh-book-2009,lubich-moldyn-2008},
financial mathematics \cite{sloan-quasi-MC-1998,wang-high-dim-finance-2006} and many others.
To work with $d$-dimensional arrays is a challenging problem due to the \emph{curse of dimensionality} \cite{beylkin-high-2005}: the number of elements of a tensor
grows exponentially with the number of dimensions $d$, and so does the complexity to work
with fully populated tensors.
For $d$ of order tens or hundreds some other approaches are needed, for example, special
low-parametric representations or \emph{formats}.
As soon as such a format comes with fast linear algebra algorithms, such as additions, Matrix-by-Vector multiplications and scalar products, any of the classical iterative methods can be implemented straightforwardly.
The first problem is that the \emph{effective} amount of unknowns, required to store vectors involved in computations in a chosen format might grow in general arbitrary during the solution process.
Second, most of classical methods are proved to be convergent in the exact arithmetics, and their behavior with approximate computations, arising from the use of formats is under the question. If the first issue depends essentially on a problem and has to be considered with a strong connection to a particular application, the second one gives more chances to be described in quite a general case. As for the GMRES method, such consideration will be presented in this paper.

Several formats have been proposed
to represent a tensor in a data-sparse way. They include canonical and Tucker formats, the two
formats with well-established properties and application areas, see the review \cite{kolda-review-2009} for more details.
They have known drawbacks. To avoid these drawbacks, the development of new tensor formats began.
In 2009 independently Hackbusch and Kuhn and later Grasedyck \cite{hk-ht-2009,gras-hsvd-2010} and Oseledets and Tyrtyshnikov \cite{ot-tt-2009} proposed two
(slightly different) hierarchical schemes for the tensor approximation, the $\mathcal{H}$-Tucker and Tree Tucker formats.
These formats depend on specially chosen \emph{dimension trees} and require recursive procedures.
To avoid the recursion, it was proposed to use a simple matrix product form of the decomposition \cite{osel-compact-2009pre,osel-tt-2011}, that was called the \emph{Tensor Train format}, or simply the TT-format.

A tensor $A$ is said to be in the TT-format, if its elements are defined by a formula
 \begin{equation}\label{solve:tt}
  A(i_1,\ldots,i_d) = G_1(i_1) \ldots G_d(i_d),
 \end{equation}
where $G_k(i_k)$ is an $r_{k-1} \times r_k$ matrix for each fixed $i_k,~1 \leq i_k \leq n_k$.
To make the matrix-by-matrix product in \eqref{solve:tt} scalar, boundary conditions $r_0 = r_d = 1$
are imposed. The numbers $r_k$ are called \emph{TT-ranks} and $G_k(i_k)$ are \emph{cores} of the TT-decomposition of a given tensor.
If $r_k \leq r, n_k \leq n$, then the storage of the TT-representation requires $\leq dnr^2$ memory cells. If $r$ is small, then this is much smaller than the storage of the full array, $n^d$.

One can go even further and introduce the \emph{Quantized TT} (QTT) format \cite{osel-2d2d-2010,khor-qtt-2011}: if the mode sizes are equal to $2^p$, we can reshape a given tensor to the $2\times 2 \times \cdots 2$ tensor with higher dimension, but all mode sizes are equal to $2$, and then apply the TT approximation to this new tensor. The storage in this case is estimated as $\mathcal{O}(d r^2 \log n)$.

The TT-format is
stable in the sense that the best approximation with bounded TT-ranks always exists and
a quasioptimal approximation can be computed by a sequence of SVDs of
auxiliary matrices \cite{osel-compact-2009pre,osel-tt-2011}.

The TT-format comes with all basic linear algebra operations.
Addition, matrix-by-vector product, elementwise multiplication can be implemented in linear $d$ and polynomial in $r$ complexity with the result also in the TT-format \cite{osel-compact-2009pre,osel-tt-2011}.
The problem is that after such operations TT-ranks grow.
For example, the TT-ranks of the sum of two tensors are equal (formally) to the sum of the TT-ranks of the addends.
In the case of the matrix-by-vector product the result is also in the TT-format with the TT-ranks of matrix and vector multiplied.
After several iterations, the TT-ranks will become too large, thus the tensor \emph{rounding},
or truncation is needed: a given tensor $A$ is approximated by another tensor $B$
with minimal possible TT-ranks with a prescribed accuracy $\varepsilon$ (or a fixed maximal rank $R$):
$$
B = \mathcal{T}_{\varepsilon,R}(A), \quad \mbox{so that}~ ||A-B||_F \leq \varepsilon ||A||_F,~\mbox{and/or}~\operatorname{rank}(B) \le R.
$$
The rounding procedure in the TT-format can be implemented in $\mathcal{O}(dnr^3)$ operations \cite{osel-compact-2009pre,osel-tt-2011}.

In this work we implement the adapted GMRES solver using the TT arithmetics and truncations. It is worth to mention the previous papers devoted to Krylov methods with tensor computations:
\cite{TobKress-Krylov-2010} (a FOM-like method in the case of a Laplace-like matrix),
\cite{TobKress-Param-2010pre} (Richardson, PCG and BiCGStab methods with application to parametric and stochastic PDEs),  and, the closest to our work, a projection method for linear systems in the $\mathcal{H}$-Tucker format which was proposed in \cite{balgras-Htuck_gmres-2010}. Our GMRES method is slightly different. First, the $\mathcal{H}$-Tucker method from \cite{balgras-Htuck_gmres-2010} uses the projectors with equal sizes:
$$
Ax=b \rightarrow W_m^{\top} A V_m y = W_m^{\top} b, \quad V_m \in \mathbb{C}^{n \times m}, ~W_m = AV_m.
$$
A very important feature of the GMRES method is the rectangular Hessenberg matrix computed via projections to subspaces with \emph{different} dimensions:
$$
\bar H_m = V_{m+1}^{\top} A V_m.
$$
(In this sense, the mentioned above method is a certain realization of \emph{geometric} minimal residual method, the linear span of (nonorthogonal) vectors $W_m$ contains all Krylov vectors from $V_{m+1}$ except the first one).
Moreover, we provide the error and convergence analysis with respect to the tensor rounding using the theory of inexact GMRES, and the relaxation strategies to reduce TT ranks and improve the performance. A convergence estimate was also provided for the tensor CG-type method for the eigenvalue problems in \cite{lebedeva-tensorcg-2010}.
A part of our paper (Section \ref{sec:precs_and_smooths}) is devoted to the role of Matrix-by-Vector and Preconditioner-by-Vector multiplications in approximate computations, showing the differences between GMRES and CG methods.

Note that the particular choice of the TT format in this paper is not important for the general theory and is due to the simplicity, convenience and robustness of the TT format in a wide class of problems.
The performance improvements considered below, arising from the inexact Krylov theory, were also successfully applied for the tensor version of GMRES in the \emph{Tucker} format by Dmitry Savostyanov. Numerical examples were presented on the Workshop on Tensor Decompositions and Applications (TDA 2010), Monopoli, Bari, Italy.

The rest of the paper is organized as follows. In the next section we briefly review the scheme of the GMRES method.
In Section \ref{sec:precs_and_smooths} we discuss the influence of preconditioners, especially in the case of errors caused by tensor roundings via SVD.
In Section \ref{sec:inexact_gmres} we explain the inexact GMRES theory, provide the error analysis for the approximate TT-GMRES and the whole algorithm.
And in the last section five we present several numerical examples, and compare also two methods: TT-GMRES and the DMRG-ALS linear solver from \cite{DoOs-dmrg-solve-2011}.

\section{Exact GMRES method in standard arithmetics}\label{sec:gmres_basic}
In this section we briefly recall the GMRES method following \cite{Saad_Schultz-gmres-1986}.
We are going to investigate, how influence the errors arising from the tensor roundings to the convergence of the methods. Moreover, if we are solving discretized PDEs, we have to consider carefully, which norm of the residual one shall use.

Let us present the basic properties of this method.
Suppose we are going to solve a linear system
$$
Ax = b, \quad A \in \mathbb{C}^{n \times n}, \quad x,b \in \mathbb{C}^{n}.
$$
The method is based on the minimization of the residual functional $||b-Ax||$ on the Krylov subspaces:
$$
\mathcal{K}_k = \left\{b, Ab, A^2b, ..., A^{k-1} b \right\}.
$$
To build the Krylov basis one uses the Arnoldi process (see the algorithm \ref{alg:GMRES_gmres}),
which is nothing else but the Gramm-Shmidt orthogonalization method applied to the Krylov vectors. After $k$
steps we have the orthonormal vectors $V_{k+1}$ and $k+1 \times k$ matrix $\bar H_k = [h_{i,j}]$.
Now we have to obtain a correction to the solution.

The vectors $v_i$ possess the following property:
\begin{equation}
A V_k = V_{k+1} \bar H_k.
\label{eqn:GMRES_reduct}
\end{equation}
Suppose the initial guess $x_0$ is given, the initial residual $r_0 = b - Ax_0$. We are to solve the least squares problem
$$
\min\limits_{z \in \mathcal{K}_k} ||f-A(x_0+z)|| = \min\limits_{z \in \mathcal{K}_k} ||r_0-Az||.
$$
Now put $z=V_k y$, reformulate the functional for the vector $y$, which is small, if $k \ll n$:
$$
J(y) = ||\beta v_1 - A V_k y||,
$$
where $\beta = ||r_0||$. Taking into account \eqref{eqn:GMRES_reduct}, we obtain
$$
J(y) = ||V_{k+1}(\beta e_1 - \bar H_k y)|| = ||\beta e_1 - \bar H_k y||,
$$
since $||V_{k+1}||=1$ due to the orthogonality, $e_1$ is the first identity vector of size $k+1$. Now we write the correction to the solution:
$$
x_k = x_0+V_k y_k, \quad y_k = \arg\min\limits_{y} ||\beta e_1 - \bar H_k y||.
$$

In the Arnoldi process the number of basis vectors grows up to the size of a matrix, resulting in a prohibitive amount of memory and computational time. To avoid this, one uses GMRES with \emph{restarts}: every $m$ steps the current solution is taken as the initial guess, and the algorithm restarts. The overall scheme of the GMRES(m) is given in the Algorithm \ref{alg:GMRES_gmres}.
\begin{algorithm}[h!]
 \caption{GMRES(m)} \label{alg:GMRES_gmres}
 \begin{algorithmic}[1]
  \REQUIRE Matrix $A$, right-hand side $b$, initial guess $x_0$, stopping tolerance $\eps$.
  \ENSURE Approximate solution $x_m:~||Ax_m-b|| \le \eps$.
  \STATE Start: compute $r_0=b-Ax_0$, $v_1=r_0/||r_0||$.
  \STATE Iterations:
     \FOR{$j=1,2,...,m$}
        \STATE $h_{i,j}=(Av_j, v_i)$, $i=1,2,...,j$, \COMMENT{Arnoldi process}
	\STATE $\hat v_{j+1} = Av_j - \sum\limits_{i=1}^j h_{i,j} v_i$,
	\STATE $h_{j+1,j} = ||\hat v_{j+1}||$, and
        \STATE $v_{j+1} = \hat v_{j+1} / h_{j+1,j}$.
     \ENDFOR
  \STATE Assemble the matrix $\bar H_m = [h_{i,j}]$, $j=1,...,m$, $i=1,...,j+1$.
  \STATE Compute the least-squares solution: $y_m = \arg\min\limits_{y} ||\beta e_1 - \bar H_m y||$,
  \STATE $x_m = x_0+V_m y_m$. \COMMENT{Update}
  \STATE Restart: compute $r_m = b-Ax_m$. Stop if $||r_m|| \le \eps$.
  \STATE Otherwise set $x_0=x_m$, $v_1=r_m/||r_m||$ and go to 2.
 \end{algorithmic}
\end{algorithm}

One of the nice properties following from \eqref{eqn:GMRES_reduct} is a cheap way to compute the residual:
$$
||\beta e_1 - \bar H_k y_k|| = ||b - A(x_0+V_k y_k)||,
$$
so we can check the stopping criteria using only small vectors $e_1, y_k$ and matrix $\bar H_k$.

For the \emph{exact} GMRES the following property is shown in \cite{Saad_Schultz-gmres-1986}:
\begin{theorem}\label{thm:gmres_exact_sol}
The solution $x_j$, obtained on $j$-th GMRES step is exact if and only if $\hat v_{j+1}=0 \Leftrightarrow h_{j+1,j}=0$.
\end{theorem}
In the following we will consider the \emph{inexact} GMRES, for which this theorem does not take place.

\section{Role of preconditioners and smoothers}\label{sec:precs_and_smooths}
A well known way to improve the convergence of an iterative method is to use a preconditioner:
$$
Ax = b \rightarrow MAx = Mb, \quad \mbox{or}~AMy=b,~x=My,
$$
which may shrink the spectrum of a matrix to a small interval (for the discretized PDE problems one usually requires a mesh-independent spectral interval), or make clusters of eigenvalues (see \cite{tee_clustering_1996,dkot-P2-2011}).

The first formula is called a \emph{left} preconditioner, the second is a \emph{right} one.
The main difference of these approaches (at least, for symmetric matrices) is what residual is computed and considered:
the first works with $||MAx - Mb||$, i.e. the \emph{preconditioned} residual, the second - with the real residual $||Ax-b||$.
Usually in the solution process the residual-based stopping criterion is used, and in the case of the left preconditioner the residual has to be computed explicitly (whereas the norm of the preconditioned residual can be computed rapidly in the GMRES method).
However, if the preconditioner is close enough to the inversed matrix ($MA = I + E$, $||E|| \ll 1$), it might be worth to use the preconditioned residual, since it provides information on the \emph{solution error}, which is more important (and relevant), than the residual. Indeed,
$$
MAx-Mb = (I+E)x - (I+E)A^{-1}b = (x - A^{-1}b) + Ex-EA^{-1}b \approx x-A^{-1}b.
$$
In some cases the norm of $E$ might be even greater than 1, but if it does not depend on a grid size,
the preconditioned residual still gives relevant information on the error,
in the sense that the constants of
equivalence $c_1 ||x-A^{-1}b|| \le ||MAx-Mb|| \le c_2 ||x-A^{-1}b||$ do not depend on $h$. Moreover, on the usual scales of errors arising in tensor arithmetical roundings ($10^{-4}-10^{-6}$),
the real residual might give no information on the convergence at all.

So, consider the tensor rounding in the following form. Suppose a tensor $u$ is given,
and consider a low-rank approximation
$$
\tilde u = u + \eps.
$$
Since the correction $\eps$ is composed from the last singular vectors of TT-blocks of $u$, it contains in general harmonics with higher frequencies, than $\tilde u$.
Let us illustrate it on the following example.
Consider a 1D function $u$ on the interval $[-1,1]$ and its Fourier decomposition:
$$
u(x) = \alpha_0 + \sum\limits_{m=1}^{\infty} \alpha_m \cos(\pi m x) + \beta_m \sin(\pi m x).
$$
Take a partial sum of this sequence as an approximation.
$$
\tilde u(x) = \alpha_0 + \sum\limits_{m=1}^{R} \alpha_m \cos(\pi m x) + \beta_m \sin(\pi m x).
$$
From the Parseval's theorem it is known, that if the coefficients are computed as follows:
$$
\alpha_m = \frac{(u,\cos(\pi m x))}{(\cos(\pi m x), \cos(\pi m x))}, \quad \beta_m = \frac{(u,\sin(\pi m x))}{(\sin(\pi m x), \sin(\pi m x))},
$$
then the approximation with harmonic functions $\tilde u$ is \emph{optimal} in the $L_2$-norm.
The approximation error is written as the following sum:
$$
u - \tilde u = \sum\limits_{m=R+1}^{\infty} \alpha_m \cos(\pi m x) + \beta_m \sin(\pi m x),
$$
i.e. it contains harmonics with the frequencies greater than $R$. The singular value decomposition provides the optimal rank-$r$ approximation to a matrix in the $2$-norm, and the same phenomenon occurs.

Consider now action of the second derivative operator $d^2/dx^2$ to the approximated function $\tilde u$. It reads
$$
\dfrac{d^2 \tilde u}{dx^2} = \dfrac{d^2 u}{dx^2} + \sum\limits_{m=R+1}^{\infty} \pi^2 m^2 \alpha_m \cos(\pi m x) + \pi^2 m^2 \beta_m \sin(\pi m x),
$$
and
$$
\left\| \dfrac{d^2 \tilde u}{dx^2} - \dfrac{d^2 u}{dx^2} \right\| \ge \pi^2 R^2 \left\|\tilde u - u \right\|.
$$
The approximation $\tilde u$ in the $L_2$ norm might provide a sufficient accuracy $\eps$, but the error in the second derivative is in the order of $R^2 \eps$, which might be prohibitively large. The discretization matrix $A$ of the second derivative operator has the norm $\mathcal{O}(1/h^2)$, so
$$
||A\tilde u - Au|| \le ||A|| \eps = \mathcal{O}\left(\frac{1}{h^2} \eps\right) \gg \eps.
$$

When we consider relative quantities, if the accuracy of the linear system solution is
$||x - A^{-1}b||/||x|| = \eps$, the residual norm might be
$||Ax-b||/||b|| = \mathcal{O}(\cond(A)~\eps)$. If the tensor rounding accuracy is $\eps=10^{-5}$,
and the problem is discretized on 1000 grid points, then $\cond(A) = \mathcal{O}(10^6)$ and
$||Ax-b||/||b|| = \mathcal{O}(10)$. The $L_2$-norm accuracy of the order $10^{-5}$ might be sufficient, but one can not control it, as the residual norm is greater than 1.

So, for tensor iterative methods, use of a preconditioner is important not only for the convergence
acceleration, but also to keep the equivalence constants between the error and the residual in the order of 1. Moreover, it is important to use the left preconditioner, when we multiply first the stiffness matrix on a vector, and then the preconditioner, which reduces the errors introduced by tensor roundings (\emph{smoother}).
In this sense, the Conjugate Gradient method is not very efficient. Indeed, consider the PCG algorithm (it works in terms of the right preconditioner $AMM^{-1}x=b$, see e.g. \cite{Hageman-Applied-iteratives-1981}):
$$
\begin{array}{l}
r_0 = b-Ax_0, \quad p_1 = M r_0 \\
\alpha_i = (r_{i-1}, M r_{i-1})/(Ap_i, p_i) \\
x_i = x_{i-1} + \alpha_i p_i \\
r_i = r_{i-1} - \alpha_i A p_i \\
\beta_i = (r_i, M r_{i})/(r_{i-1}, Mr_{i-1}) \\
p_{i+1} = M r_i + \beta_i p_i.
\end{array}
$$
In the formulas for $\alpha_i$ and $r_i$ the last operation before the scalar product or linear combination is the MatVec with the stiffness matrix $A p_i$, which might be computed with a very large error, if $p_i$ is rounded in the $L_2$ norm. As a result, new iterands are computed with such a large error, and the method diverges. In this sense, it is much more efficient to use the GMRES method with the left preconditioner, when the next Krylov vector is computed as $v_{i+1} = MAv_i$. Numerical experiments conducted show, that the roundings in this operation can be even applied sequentially:
$$
v_{i+1} = \mathcal{T}_{\eps,R}(M \mathcal{T}_{\eps,R}(Av_i))
$$
without corrupting the final result. It is especially important if the preconditioner is combined from several matrices, e.g. the preconditioner from \cite{dkot-P2-2011}, when the whole matrix $MA$ can not be computed explicitly due to very high ranks, and multiplications are applied during several successive implicit procedures, which provide approximate products (for example, the DMRG-based TT-MatVec, see \cite{khos-dmrg-2010,DoOs-dmrg-solve-2011} for the DMRG schemes).

\section{Relaxation strategies, inexact GMRES, and the TT-GMRES algorithm}\label{sec:inexact_gmres}
The inexact Krylov methods theory \cite{Giraud-gmres_note-2004,Simoncini-Theory-inexact-Krylov-2003,Eshof-Inexact-Krylov-Subspace-2004} allows us to estimate influence of the noise arising from the tensor
roundings, on the GMRES convergence. Moreover, the performance of tensor methods essentially depends on the TT ranks of the intermediate vectors.
It turns out in practice, that if we truncate all the vectors with the same accuracy, the ranks of the last Krylov vectors, being added to the Krylov basis, increase from iteration to iteration.
Relaxation strategies,
proposed in the papers \cite{Giraud-gmres_note-2004,Simoncini-Theory-inexact-Krylov-2003,Eshof-Inexact-Krylov-Subspace-2004}
allow us to truncate the latter Krylov vectors less accurately, than the former ones, thus keeping the ranks at almost constant values, or even reducing them on last iterations.

We can consider roundings as the usual MatVecs with the perturbed matrix:
$$
\mathcal{T}_{\eps,R}(Av) = \tilde Av = (A+E) v,
$$
where $E$ is the error matrix, which, generally speaking, changes each time, when the Matrix-by-Vector multiplication is computed, and moreover, might depend on $v$. Our goal is to estimate allowed values for $||E||$ which provide the convergence until some desired accuracy.
It can be proved, that during the Krylov iterations the norm of the error matrix $||E||$ can be increased, i.e. the vector operations on the last iterations can be computed with worse accuracy.

A linear system reduction to the Krylov subspace \eqref{eqn:GMRES_reduct} for the inexact GMRES is written as:
\begin{equation}
\begin{bmatrix}(A+E_1) v_1 & \cdots & (A+E_m) v_m \end{bmatrix} = V_{m+1} \bar H_m.
\label{eqn:GMRES_reduct_inexact}
\end{equation}
Notice that $V_m$ is no more a basis in the exact subspace $\mathcal{K}_m$. The algorithm remains the same, as for the exact GMRES, but the minimization of $||\beta e_1 - \bar H_m y_m||$ does not lead any more to the minimization of the real residual $||b - Ax_m||$. From \eqref{eqn:GMRES_reduct_inexact} follows, that we are solving the following optimization problem:
$$
\min\limits_{q \in R(W_m)} ||r_0-q||, \quad r_0 = b - (A+E_0)x_0,
$$
and $R(W_m)$ is a linear span of the vectors $W_m = V_{m+1} \bar H_m$.
So, in fact we are minimizing the \emph{computed} approximate residual
$$
||\tilde r_m|| = ||r_0-q_m|| = |h_{m+1,m} e_m^{\top} y_m|.
$$
If the new Krylov vector $(A+E_m) v_m$ on some iteration appears to be almost linearly dependent with the basis $V_m$, the quantity $h_{m+1,m}$ is small, and the approximate residual $||\tilde r_m||$ is also small,
which we can consider as a convergence of the method. But the real residual $||r_m||$ might be significantly
larger, and the estimate of the difference $||r_m-\tilde r_m||$ will be considered below.

We formulate the main theorem which comes from \cite{Simoncini-Theory-inexact-Krylov-2003}:
\begin{theorem}\label{thm:GMRES_relax_strat}
Suppose some $\eps>0$ is given, for the system $Az = r_0$ $m$ GMRES iterations are conducted, the computed residual $\tilde r_m$ was obtained. Then if for any $i \le m$ holds that
$$
||E_i|| \le \dfrac{\sigma_m(\bar H_m)}{m} \dfrac{1}{||\tilde r_{i-1}||} \eps,
$$
where $\sigma_m(\bar H_m)$ is a minimal singular value of the Hessenberg matrix of the reduced system, for the real residual $r_m = r_0 - A z_m$ the following estimate holds:
$$||r_m - \tilde r_m|| \le \eps.$$
\end{theorem}
We refer for the proof to \cite{Simoncini-Theory-inexact-Krylov-2003}.

So, the accuracy of the MatVec computation can be relaxed inversely proportional to the current residual,
and if the process is stopped (in the case of stagnation, or if the computed residual becomes smaller than the stopping tolerance), the real residual will differ from the computed one on the quantity not greater in the norm than $\eps$, i.e. the convergence of the method is controlled.

In order to obtain scale-independent estimates (i.e. the same for the systems $Ax=b$ and $\alpha Ax = \alpha b$), one usually consider the \emph{relative} residual, and the corresponding stopping criteria, for example,
$$
\dfrac{||r_i||}{||b||} \le \epsilon.
$$
In the same way one can consider the difference between the real and computed residuals: $||r_m - \tilde r_m||/||b|| \le \epsilon$.
In this case the result of Theorem \ref{thm:GMRES_relax_strat} can be reformulated for the relative quantities, taking into account that $\eps = \epsilon ||b||$:
$$
\dfrac{||E_i||}{||A||} \le \dfrac{\sigma_m(\bar H_m)}{m ||A||} \dfrac{1}{||\tilde r_{i-1}||/||b||} \epsilon,
$$

The minimal singular value of $\bar H_m$ can be estimated from the minimal singular value of $A$:
$$
\sigma_m (\bar H_m) \ge \sigma_n(A) - \left\|\begin{bmatrix} E_1 v_1 & \cdots & E_m v_m \end{bmatrix}\right\|,
$$
so we formulate the following relaxation strategy for the MatVec error:
\begin{corollary}
Suppose $m$ GMRES iterations are conducted. If for any $i \le m$ the relative error introduced in the Matrix-by-Vector multiplication is bounded by the following rule:
\begin{equation}
\dfrac{||E_i||}{||A||} \le \dfrac{1}{m \cond(A)} \dfrac{1}{||\tilde r_{i-1}||/||b||} \epsilon,
\label{eqn:relative_relax}
\end{equation}
than the real relative residual and the computed one are connected with $\frac{||r_{m}||}{||b||} \le \frac{||\tilde r_{m}||}{||b||} + \epsilon$.
\end{corollary}
With a good spectrally equivalent preconditioner $\cond(A)=\mathcal{O}(1)$ (notice that the matrix $A$ is considered to be already left-preconditioned here), and $m = \mathcal{O}(1)$,
in this case we can consider \eqref{eqn:relative_relax} in the following form: if
$$
\dfrac{||E_i||}{||A||} \le \dfrac{1}{||\tilde r_{i-1}||/||b||} \epsilon,
$$
then the inexact GMRES will converge to the relative residual not greater than $$m \cond(A) \epsilon.$$
This approach will be used in the numerical experiments below.

Let us write the final algorithm \ref{alg:GMRES_TT-gmres} of the tensor GMRES with relaxations.
Notice also, that in the Arnoldi process we used the modified Gramm-Shmidt algorithm, which is more stable
in the presence of the rounding errors. In addition, as the left preconditioner is used, we do not write it
explicitly, but assume that the matrix $A$ and the right-hand side $b$ are already preconditioned.

One additional thing which is important to note, is when to perform tensor rounding, either after adding all the summands in the orthogonalization and correction steps, or after each addition.
Formally, one can introduce the error only to the MatVec itself, but the orthogonality of $V_m$ must be kept despite the perturbations in \eqref{eqn:GMRES_reduct_inexact}.
Moreover, significantly different in magnitude vectors $y_j(i) v_i$ are also better to sum exactly.
The obvious drawback is the rank overhead which can be $m$ times larger than in the case of step-by-step truncations.
So when possible (small number of iterations) it is worth to perform only final truncation when the summation is ready (in the case of small mode sizes (Quantized TT) it can be easily done by the DMRG truncation (see next section) instead of the direct one from \cite{ot-tt-2009}).
\begin{algorithm}[h!]
 \caption{Relaxed TT-GMRES(m)} \label{alg:GMRES_TT-gmres}
 \begin{algorithmic}[1]
  \REQUIRE Right-hand side $b$, initial vector $x_0$ in the TT format, matrix $A$ as a tensor MatVec
  procedure $y=\mathcal{T}_{\eps,R}(Ax)$, accuracy $\eps$ and/or maximal TT rank $R$.
  \ENSURE Approximate solution $x_j:~||Ax_j-b||/||b|| \le \eps$.

  \STATE Start: compute $r_0=\mathcal{T}_{\eps,R}(b-Ax_0)$, $\beta = ||r_0||$ $v_1=r_0/\beta$.
  \STATE Iterations:
     \FOR{$j=1,2,...,m$}
        \STATE Compute the relaxed accuracy $\delta = \dfrac{\eps}{||\tilde r_{j-1}||/\beta}$.
        \STATE $w = \mathcal{T}_{\delta,R}(A v_j)$ - new Krylov vector.
	\FOR{$i=1,2,...,j$} 
	  \STATE $h_{i,j}=(w, v_i)$,
	  \STATE $w = w - h_{i,j} v_i$, \COMMENT{orthogonalization}
        \ENDFOR
	\STATE $w = \mathcal{T}_{\delta,R}(w)$. \COMMENT{compression}
	\STATE $h_{j+1,j} = ||w||$, $v_{j+1} = w / h_{j+1,j}$.
	\STATE Assemble matrix $\bar H_j = [h_{i,k}]$, $k=1,...,j$, $i=1,...,j+1$.
	\STATE Compute a solution of the reduced system: $y_j = \arg\min\limits_{y} ||\beta e_1 - \bar H_j y||$.
        \STATE Check the residual $||\tilde r_j|| = ||\beta e_1 - \bar H_j y_j||$: if $||\tilde r_j||/||b|| \le \eps$, then break.
     \ENDFOR
  \STATE Update the solution: initialize $x_{j} = x_0$,
  \FOR{$i=1,2,...,j$}
    \STATE $x_{j} = x_j + y_j(i) v_i$ \COMMENT{correction}
  \ENDFOR
  \STATE $x_{j} = \mathcal{T}_{\eps,R}(x)$ \COMMENT{compression}
  \STATE Restart: if $||\tilde r_j||/||b|| > \eps$, then set $x_0=x_j$, go to 1.
 \end{algorithmic}
\end{algorithm}

\section{Fast and accurate TT arithmetics in high dimensions}\label{sec:dmrg_round}
One class of interesting high-dimensional problems is the multiparametric problems arising in the discretized Karhunen-Loeve model for the PDEs with stochastic data. In such problem, the number of parameters is usually in the order of tens, and after the tensorisation (Quantisation), the number of dimensions is in the order of hundreds. Even for 1D physical problem, the QTT ranks scale usually linear with the number of parameters, thus keep the values 50-100. In this case, the multiply-and-compress strategy fails, because of the prohibitive complexity $\mathcal{O}(d n r^6)$.
A better alternative is to use direct minimization methods, based on the alternating directions approach, the ALS and DMRG (also known as MALS) schemes. There are several papers on the linear- and eigenvalue solvers using the DMRG scheme \cite{khos-dmrg-2010,holtz-ALS-DMRG-2010,DoOs-dmrg-solve-2011}. The simple approximation problem is discussed in these articles as well, and now there is the new one \cite{Os-mvk2-2011}, concerning specially the approximate Matrix-by-Vector product.

Unfortunately, the main disadvantage of all presented TT-DMRG methods is the tendency to underestimate ranks in essentially high-dimensional problems. Recall briefly the main sketch of the approximation via the DMRG (MALS) scheme:
\begin{enumerate}
 \item Suppose a functional $J(x)$ to minimize is given (e.g. $J(x)=||x - y||^2$).
 \item Initial guess for $x$ in the TT format is given: $x = X_1(i_1) \cdots X_d(i_d)$.
 \item Choose two neighboring cores and convolve a \emph{supercore}: $X_k(i_k) X_{k+1}(i_{k+1}) \rightarrow W_k(i_k,i_{k+1})$.
 \item Solve the reduced optimization problem for the elements of $W_k$: $\hat{W}_k = \arg\min\limits_{W_k} J(x)$.
 \item Recover the TT structure (e.g. via SVD): $\hat{W}_k \approx \hat{X}_k(i_k) \hat{X}_{k+1}(i_{k+1})$.
 \item Consider the next pair of cores, and so on..
\end{enumerate}
The rank is determined adaptively on the step 5. The ranks are not known in general, and we usually start from a low-rank initial guess, subsequently increasing them during the DMRG iterations. The problem is that if we are using the fixed $\eps$-truncation of singular values, the ranks determined become underestimated, as the dimension increases. There are two factors. First, the worst-case error accumulation in the whole tensor is $d \eps$, if the local errors in each block are bounded by $\eps$ \cite{ot-tt-2009}. Second, instead of direct compression routine from \cite{ot-tt-2009}, where the fixed cores are cores of the initial tensor, here we are working with a \emph{projection} to some tensor with blocks, which are far from the good approximation (on early iterations), and moreover, have insufficient ranks. To get rid of this, in this work we used the algorithms modified as follows:
\begin{itemize}
 \item First, set the accuracy for the local truncation to $\eps_{loc} = \eps/d$.
 \item Second, after the rank is truncated according to $\eps_{loc}$, artificially add more singular vectors (thus obtaining the truncation with increased accuracy and rank). This additional rank can even be determined adaptively, depending on the convergence of the current supercore, by comparison with the approximation from the previous iteration.
\end{itemize}
The approximation computed this way might have overestimated ranks. To reduce them to proper values, it is sufficient to conduct the last iteration with the standard truncation without including additional singular vectors (in fact, it performs like the direct compression routine, as the proper approximation is already achieved on this step, but the complexity is now $\mathcal{O}((r+r_{add})^3)$ instead of $\mathcal{O}(r^6)$, and the additional rank is usually significantly smaller than $r$).

For the DMRG-solve routine, we will show the role of the increased-rank truncation in the next section. But for the approximations and MatVecs in the TT-GMRES, we always keep it on.

\section{Numerical experiments}\label{sec:numerics}
The TT-GMRES method and the numerical experiments were implemented using the routines from TT Toolbox 2.1
(\texttt{http://spring.inm.ras.ru/osel/}) in the MATLAB R2009b and conducted on a Linux x86-64 machine with Intel Xeon 2.00GHz CPU in the sequential mode.

\begin{figure}[h!]
 \centering
 \includegraphics[width=8cm]{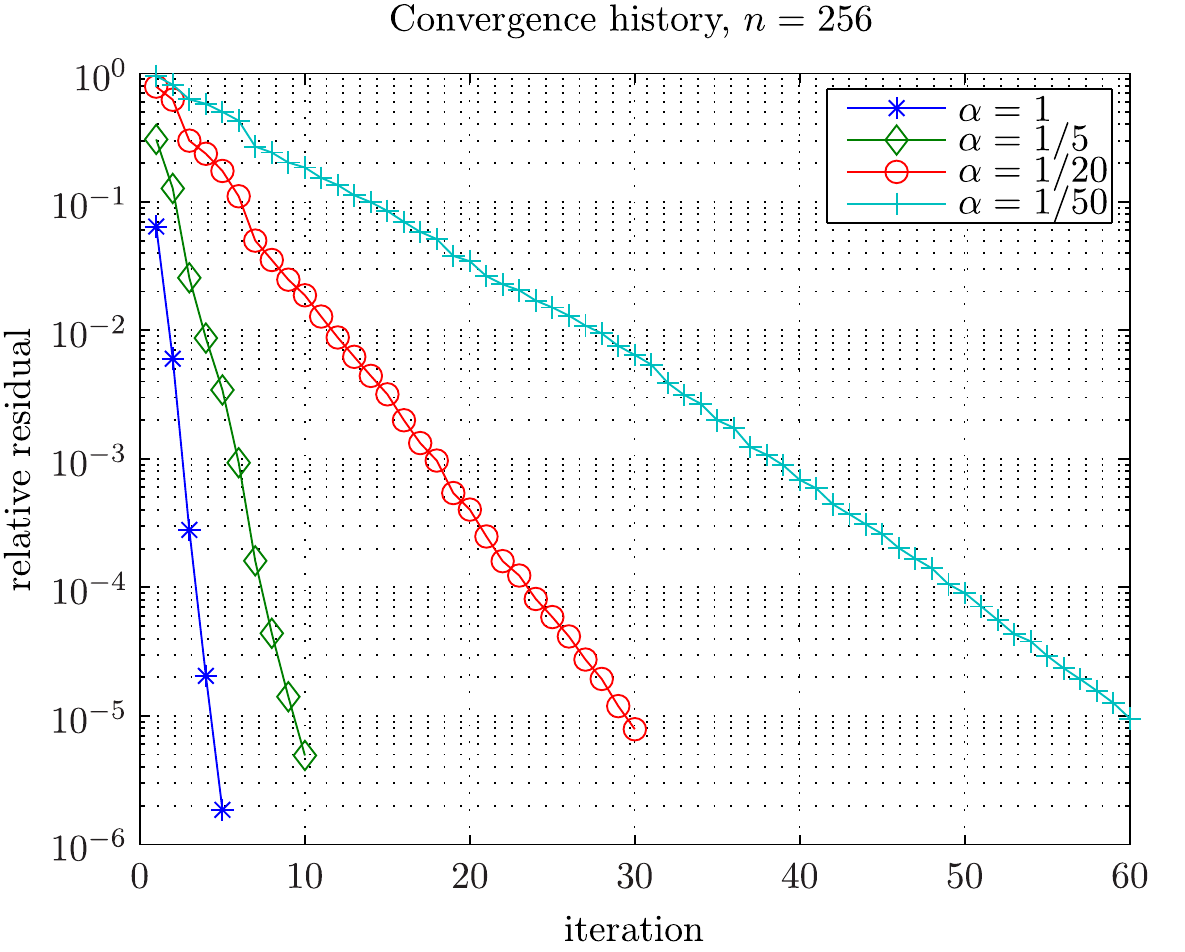}
 \caption{Convergence history for the convection example}
 \label{fig:conv:convhist}
\end{figure}

\begin{figure}[h!]
\begin{minipage}[b]{0.49\linewidth}
\centering
  \includegraphics[scale=0.52]{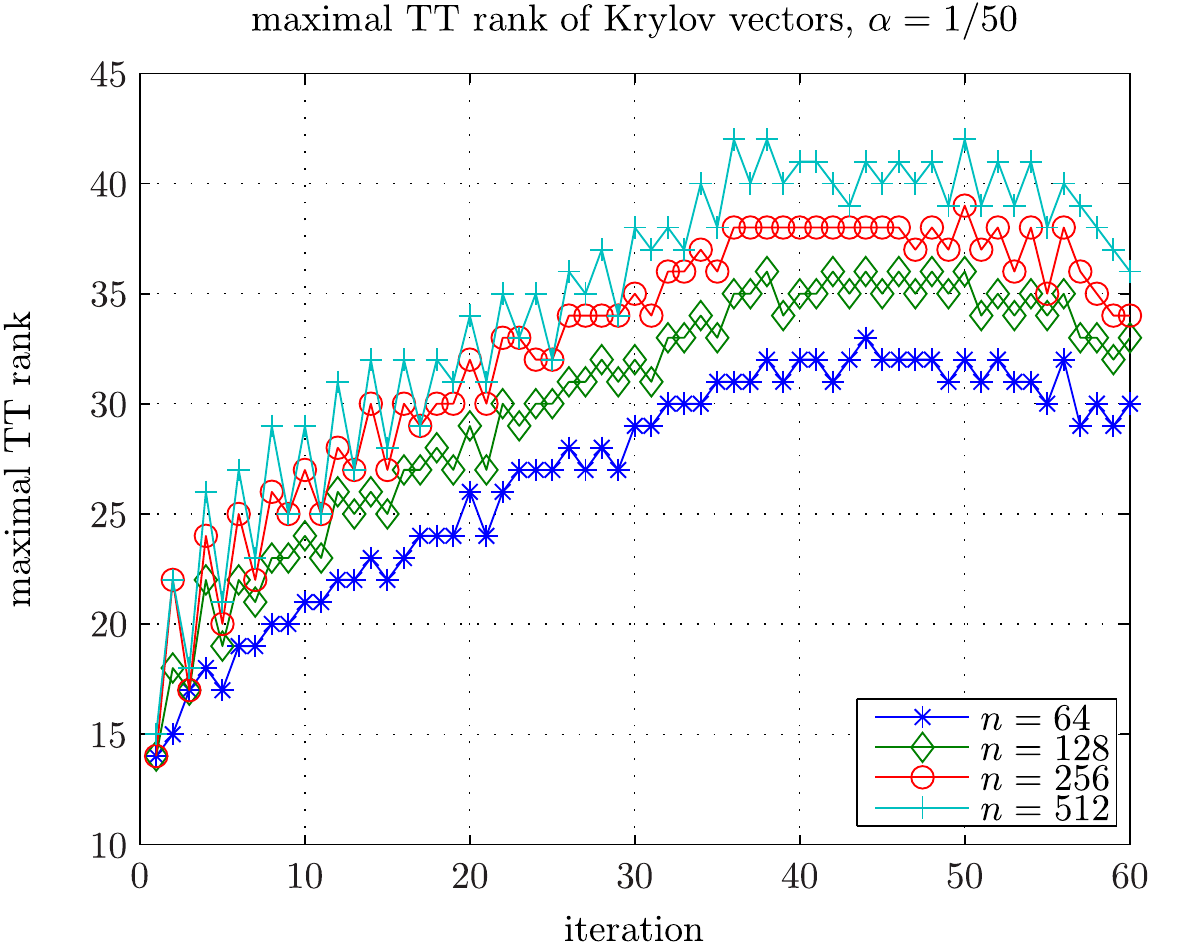}
 \caption{Maximal TT rank of the last Krylov vector, convection example}
 \label{fig:conv:rw}
\end{minipage}
\hfill
\begin{minipage}[b]{0.49\linewidth}
\centering
  \includegraphics[scale=0.52]{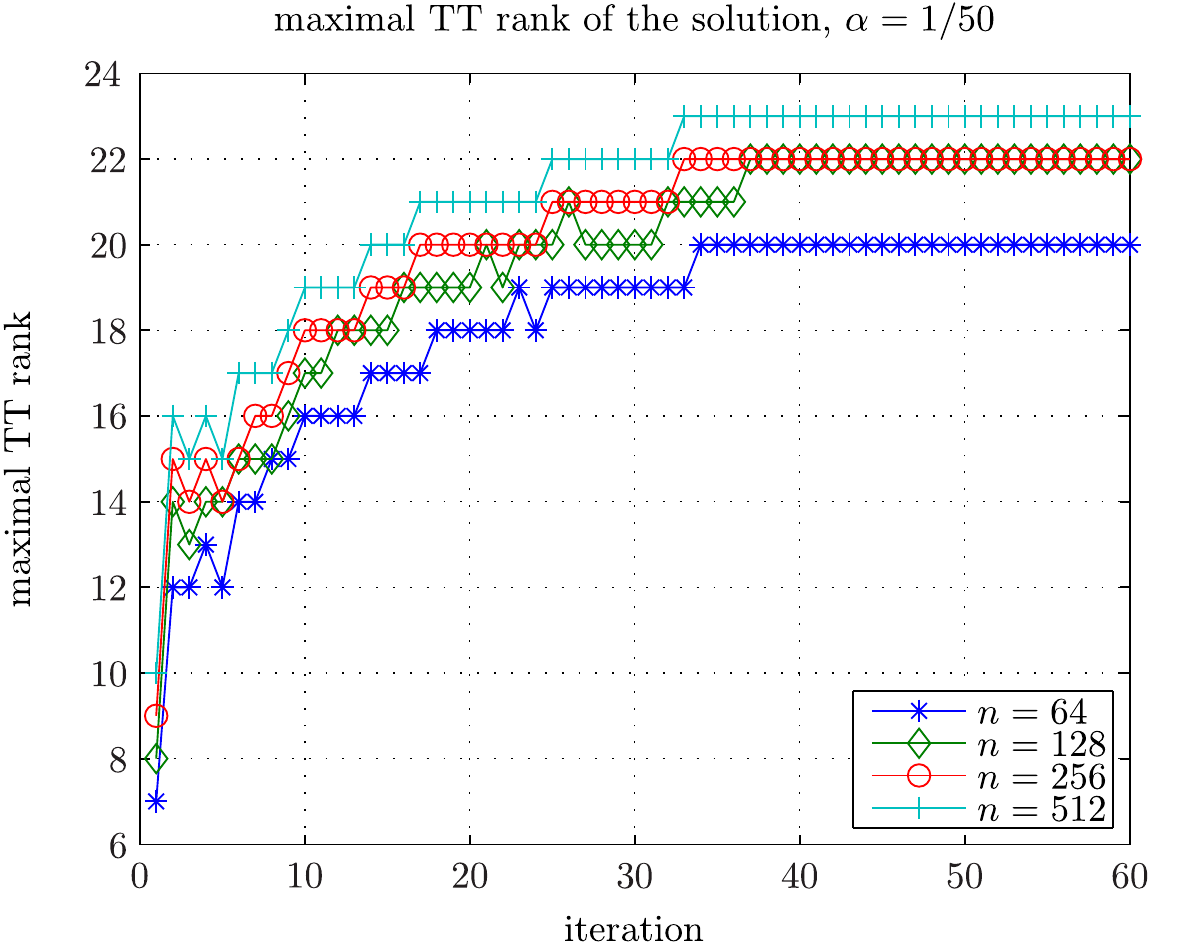}
 \caption{Maximal TT rank of the solution, convection example}
 \label{fig:conv:rx}
\end{minipage}
\end{figure}

\begin{figure}[h!]
\begin{minipage}[b]{0.49\linewidth}
\centering
  \includegraphics[scale=0.52]{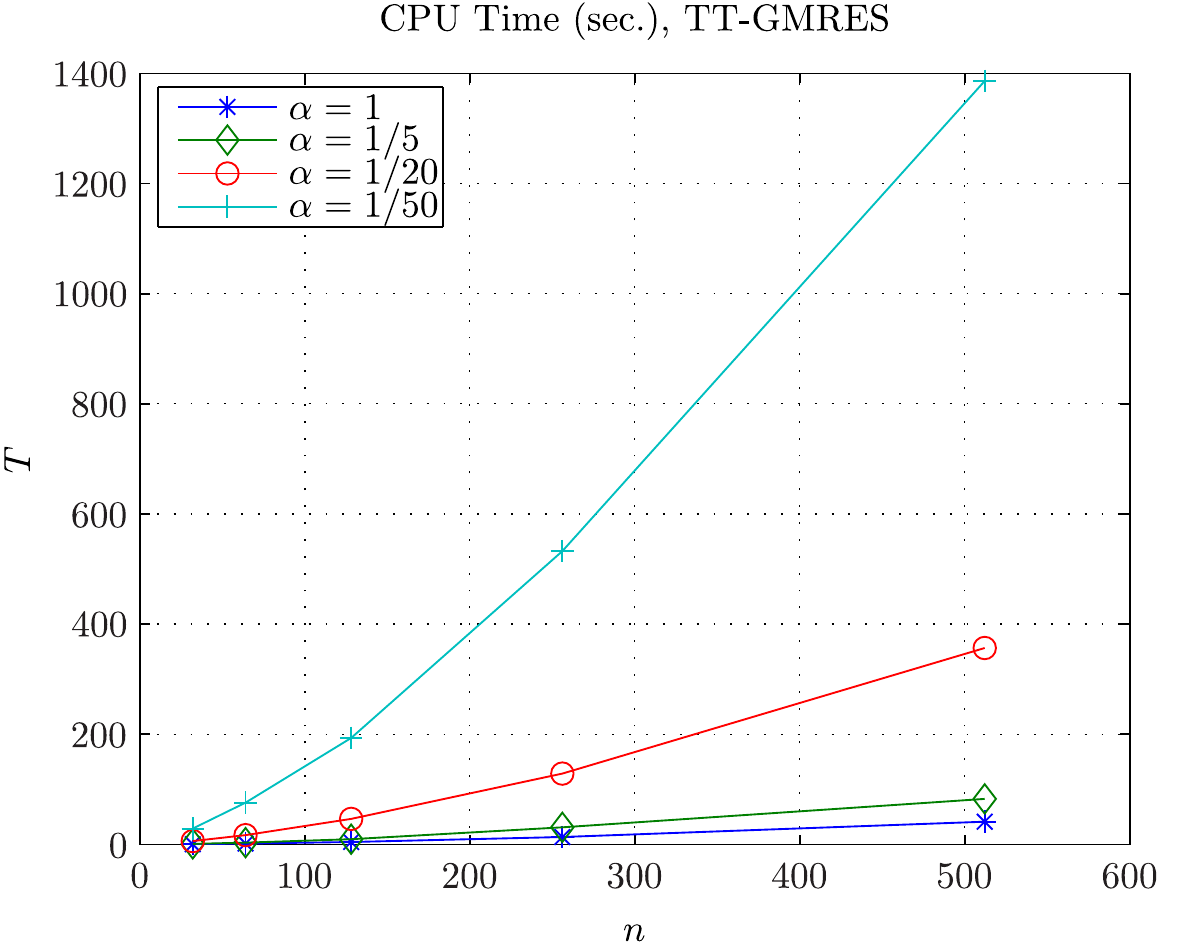}
 \caption{CPU time (sec.) of the TT solver, convection example}
 \label{fig:conv:ttts}
\end{minipage}
\hfill
\begin{minipage}[b]{0.49\linewidth}
\centering
  \includegraphics[scale=0.52]{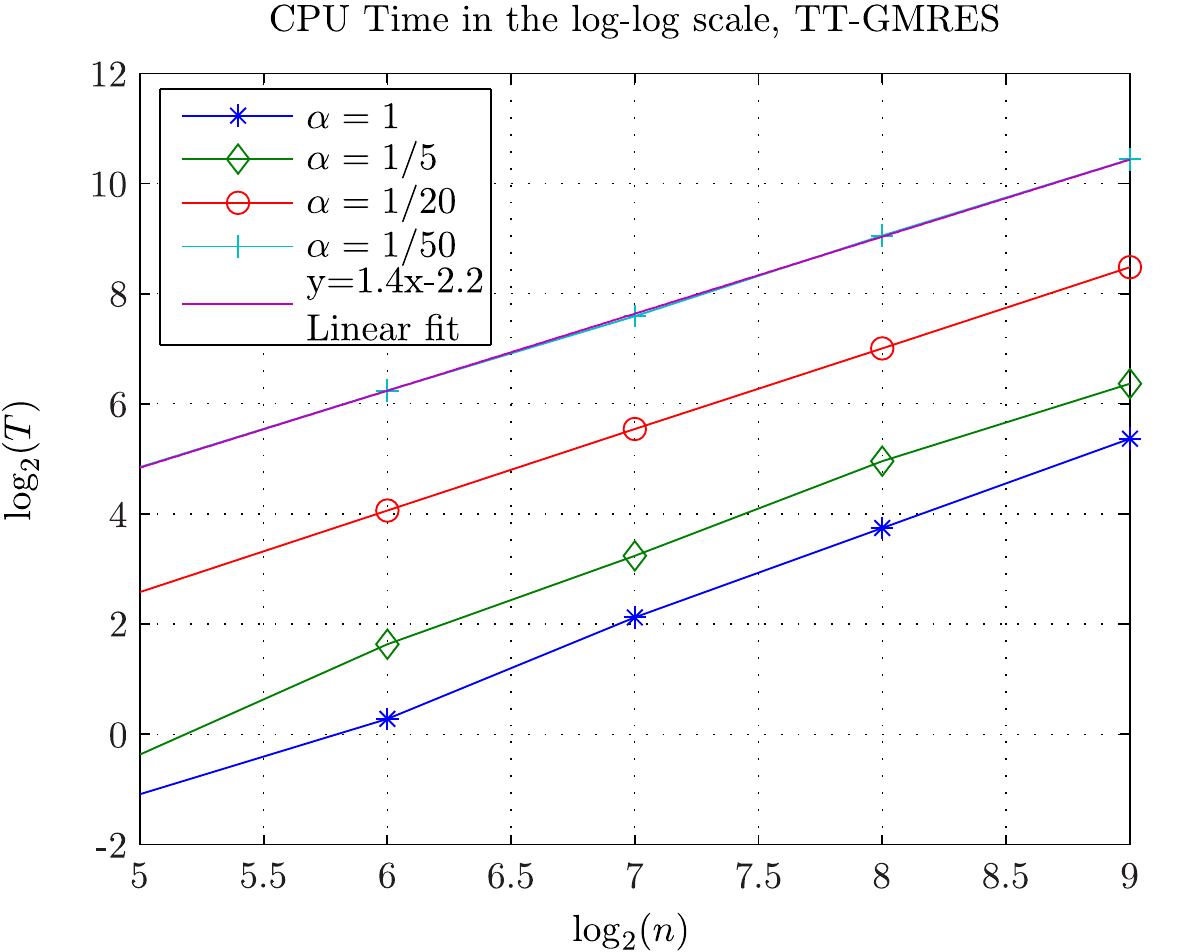}
 \caption{CPU time (sec.) of the TT solver in the log-log scale, convection example}
 \label{fig:conv:logttts}
\end{minipage}
\end{figure}

\begin{figure}[h!]
\begin{minipage}[b]{0.49\linewidth}
\centering
  \includegraphics[scale=0.52]{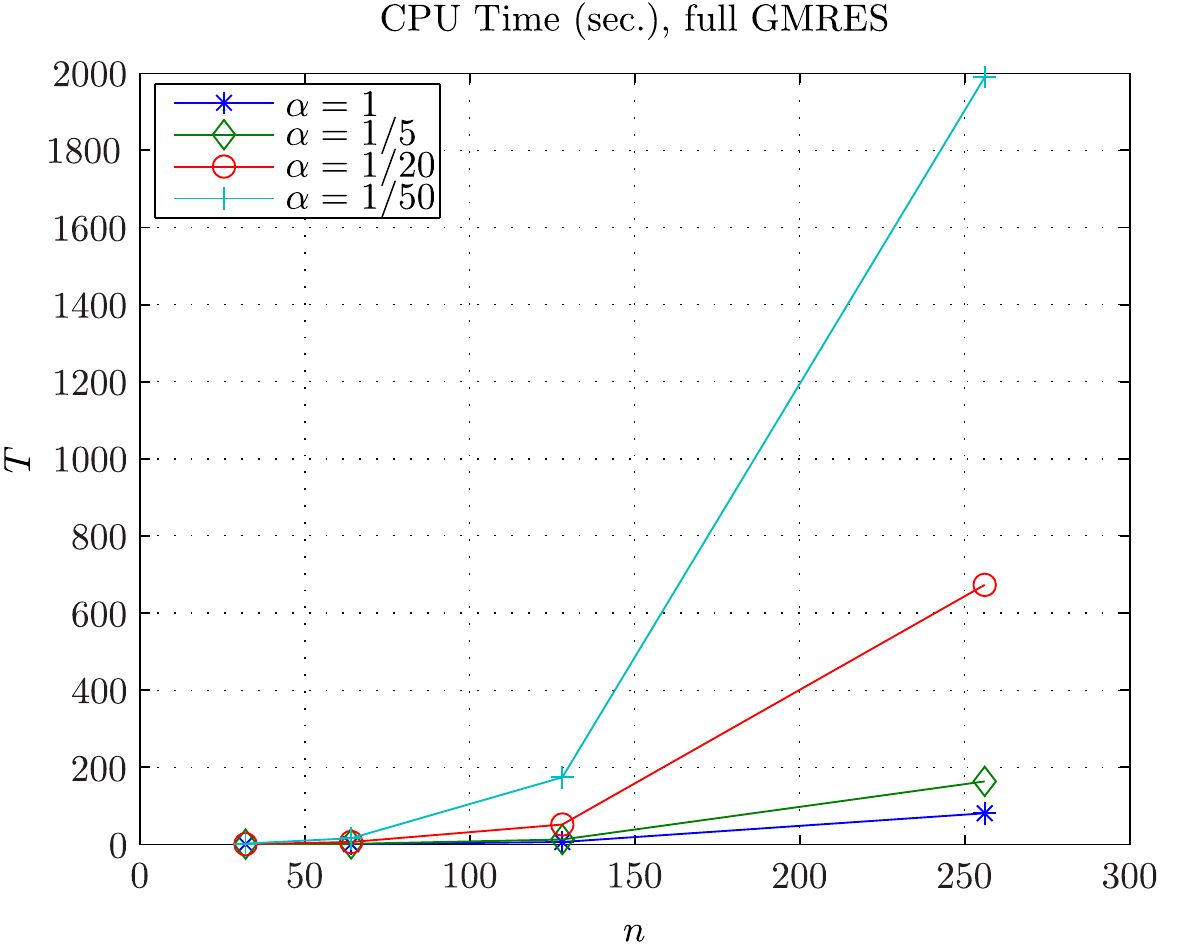}
 \caption{CPU time (sec.) of the full solver, convection example}
 \label{fig:conv:tfulls}
\end{minipage}
\hfill
\begin{minipage}[b]{0.49\linewidth}
\centering
  \includegraphics[scale=0.52]{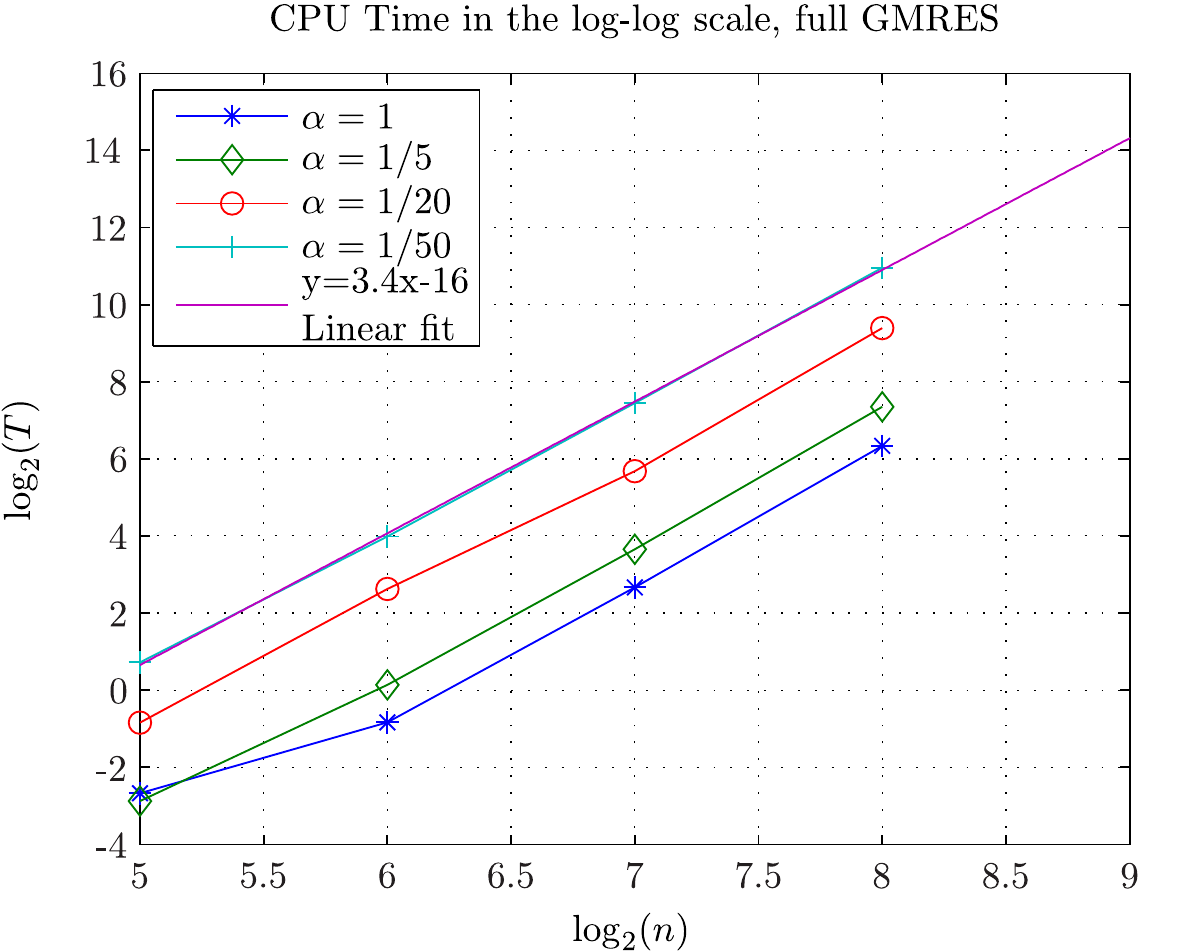}
 \caption{CPU time (sec.) of the full solver in the log-log scale, convection example}
 \label{fig:conv:logtfulls}
\end{minipage}
\end{figure}

\subsection{Convection-diffusion (Table \ref{tab:conv:it_n_alpha}, Figures \ref{fig:conv:convhist} - \ref{fig:conv:logtfulls})}
The first example is a 3D diffusion-convection problem with the recirculating wind
$$
\left\{\begin{array}{ll}
        -\alpha\Delta u + 2y(1-x^2) \dfrac{\partial u}{\partial x} - 2x(1-y^2) \dfrac{\partial u}{\partial y} = 0 & \mbox{in}~\Omega = [-1,1]^3,\\
       u_{y=1} = 1, \quad u_{\partial\Omega \backslash\{y=1\}}=0
       \end{array}\right.
$$
discretized using the central-point finite difference scheme:
$$
-\Delta ~\rightarrow~ -\Delta_h = (-\Delta^1_h) \otimes I \otimes I + I \otimes (-\Delta^1_h) \otimes I + I \otimes I \otimes (-\Delta^1_h),
$$
$$
\dfrac{\partial u}{\partial x} ~\rightarrow~ \nabla_h^x = \nabla^1_h \otimes I \otimes I, \quad \dfrac{\partial u}{\partial y} ~\rightarrow~ \nabla_h^y = I \otimes \nabla^1_h \otimes I,
$$
$$
-\Delta^1_h = \dfrac{1}{h^2}\begin{bmatrix}
                     2 & -1 & 0 \\
		     -1& 2 & -1 & 0 \\
		     0 & -1 & 2 & -1 & 0 \\
		     & & \ddots & \ddots & \ddots \\
		     &        &        &  -1 & 2
                    \end{bmatrix}, \quad
\nabla^1_h = \dfrac{1}{h} \begin{bmatrix}
                     0 & 0.5 & 0 \\
		     -0.5& 0 & 0.5 & 0 \\
		     0 & -0.5 & 0 & 0.5 & 0 \\
		     & & \ddots & \ddots & \ddots \\
		     &        &        &  -0.5 & 0
             \end{bmatrix},
$$
$h = 1/(n+1)$ is a grid size.
The scalar parameter $\alpha$ (diffusion scale) varies from $1$ to $1/50$ in the numerical tests below.

We use the TT data representation (without the QTT structure), so the TT ranks of the stiffness matrix
$$
-\alpha\Delta_h + \left(\diag\left(1-x^2\right) \otimes \diag\left(2y\right) \otimes I\right) \cdot \nabla_h^x + \left(\diag\left(-2x\right) \otimes \diag\left(1-y^2\right) \otimes I\right) \cdot \nabla_h^y
$$
are bounded by 4 (the ranks of $-\Delta_h$ are all equal to 2, see \cite{khkaz-lap-2010}).

To solve this problem efficiently, we use the inversed discrete Laplacian $-\Delta_h^{-1}$ as a preconditioner (although this is not the optimal preconditioner, and the convergence depends significantly on $\alpha$, the problem is tractable within our range of Reynolds numbers).
To implement the inversed Laplacian in the TT format we used the quadrature from \cite{khor-low-rank-kron-P1-2006,khor-low-rank-kron-P2-2006}: if
$$
\Delta_h = \Delta^1_h \otimes I \otimes \cdots \otimes I + \cdots + I \otimes \cdots \otimes I \otimes \Delta_h^{1},
$$
then
$$
\Delta_h^{-1} \approx \sum\limits_{k=-M}^M c_k \bigotimes\limits_{p=1}^d \exp(-t_k \Delta_h^{1}),
$$
where $t_k=e^{k\eta}$, $c_k=\eta t_k$, $\eta=\pi/\sqrt{M}$, with the accuracy $\mathcal{O}(e^{-\pi \sqrt{M}})$, so that $r_{\Delta^{-1}} = \mathcal{O}(\log^2(1/\eps))$.
In practice this formula can be accelerated (giving the complexity $\mathcal{O}(n~\log n)$) by using Fast Trigonometric transforms (in our case of Dirichlet boundary conditions the appropriate transform is DST-I) with all TT ranks equal to 1 \cite{dks-ttfft-2012}, and compressing only the diagonal matrix with inversed eigenvalues.

The timings of the TT solver are compared with ones of the standard full-vector GMRES solver, with the same preconditioner implemented in the full format using the trigonometric transforms as well, with the complexity $\mathcal{O}(n^3~\log n)$.

The tensor rounding accuracy for the solution is fixed to $\eps=10^{-5}$, and the accuracy for the Krylov vectors is determined according to the relaxation strategies.

First, we check the convergence properties of the preconditioner (Table \ref{tab:conv:it_n_alpha}, Fig. \ref{fig:conv:convhist}).
\begin{table}[h!]
\centering
\caption{Number of iterations versus the grid size ($n$) and diffusion scale ($\alpha$)}
\label{tab:conv:it_n_alpha}
\begin{tabular}{|c|c|c|c|c|c|c|} \hline
\backslashbox[-20pt][l]{$n$}{$\alpha$} & 1 & 1/2 & 1/5 & 1/10 & 1/20 & 1/50 \\ \hline
64					&5 & 6   & 10  & 17  & 30  & 60 \\ \hline
256					&5 & 6   & 10  & 17  & 30  & 60 \\ \hline
\end{tabular}
\end{table}
The number of iterations is stable with respect to the grid size, but grows approximately linearly with the Reynolds number.
The convergence histories for different $\alpha$ and $n=256$ are given on Fig. \ref{fig:conv:convhist}.

The behavior of the TT ranks during the iterations (we measure here the highest rank $\max\limits_{i=1,..,d-1} r_i$) of the Krylov vectors and the solution is presented on Fig. \ref{fig:conv:rw}, \ref{fig:conv:rx}, respectively.
The solution rank grows from 1 (zero tensor) to its stable value with a weak (approx. logarithmic) dependence on the grid size.
The Krylov vector rank has its maximum at the middle iterations on the finer grids (it is also important, that it grows slightly with the grid size, it will be reflected in the computational time), but near the end of the process, it begins to decrease due to the relaxed accuracy.

Now, consider the computational time of the TT-GMRES solver and the standard full GMRES method in MATLAB with the same Fourier-based preconditioner.
The CPU time of the TT solver is presented on Fig \ref{fig:conv:ttts}, and the log-log scale plot is on Fig. \ref{fig:conv:logttts}
The linear fitting on the log-log plot gives the experimental complexity rate $n^{1.4}$. The overhead with respect to the true linear complexity appears from the additional logarithmic terms in the Fourier transforms and approximately logarithmic grow of the TT ranks of the Krylov vectors, see Fig. \ref{fig:conv:rw}.

The full solver manifests the complexity rate $n^{3.4}$, which lies in a correspondence with its theoretical estimate $\mathcal{O}(n^3~\log n)$ (see Fig. \ref{fig:conv:tfulls} for the CPU time itself, and \ref{fig:conv:logtfulls} for the log-log scale).
Notice that the full solver timings are presented only for grid sizes not larger than $256$. We were not able to perform the calculations on the grid $512^3$ due to insufficient memory resources. Nevertheless, the extrapolation via the linear fit from the Fig. \ref{fig:conv:logtfulls} gives an estimate $2^{14} \sim 20000$ seconds for that experiment, which is about 15 times larger than the corresponding times of the TT solver.

\begin{figure}[h!]
\begin{minipage}[b]{0.49\linewidth}
\centering
  \includegraphics[scale=0.52]{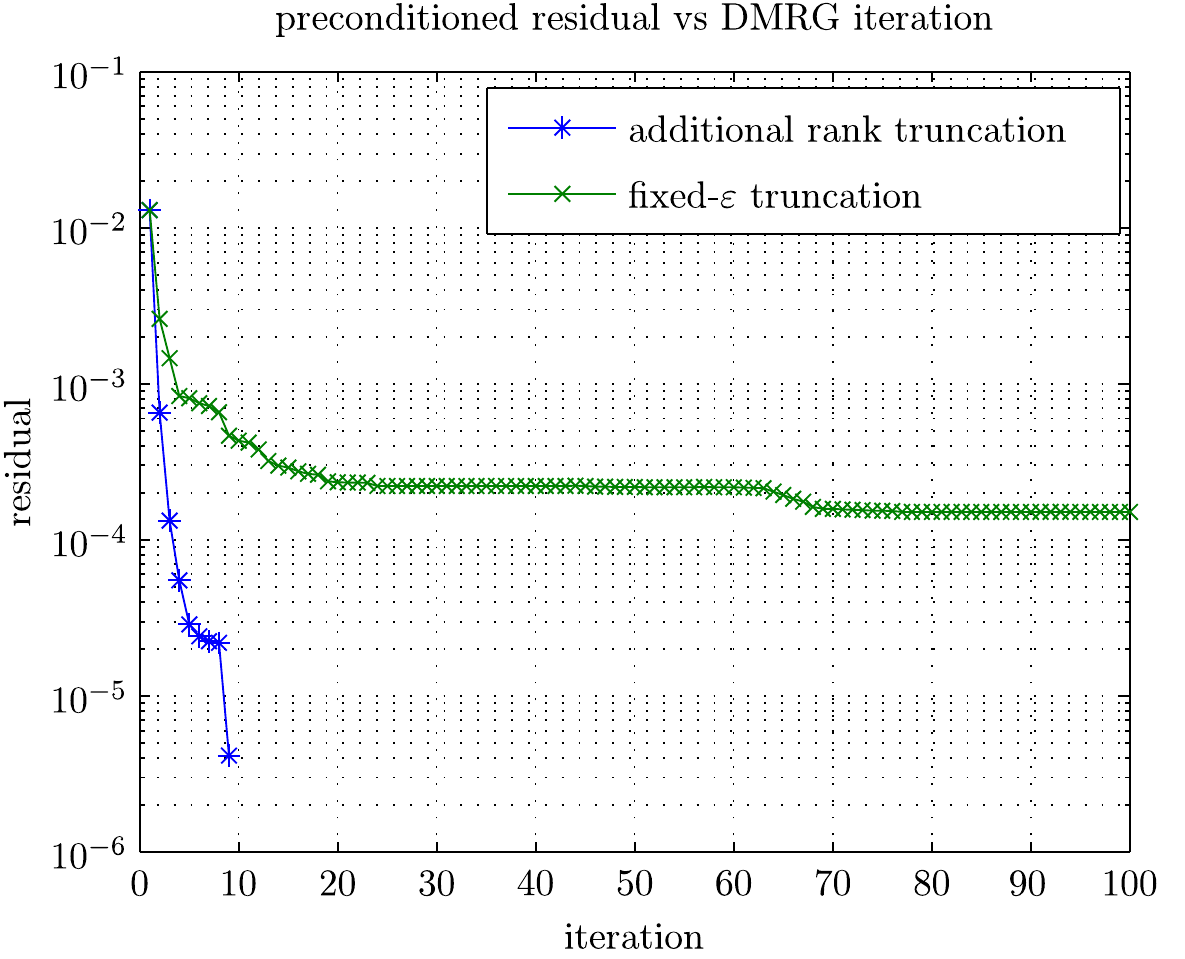}
 \caption{Relative residuals, sPDE, DMRG solvers}
 \label{fig:resid_dmrg}
\end{minipage}
\hfill
\begin{minipage}[b]{0.49\linewidth}
\centering
  \includegraphics[scale=0.52]{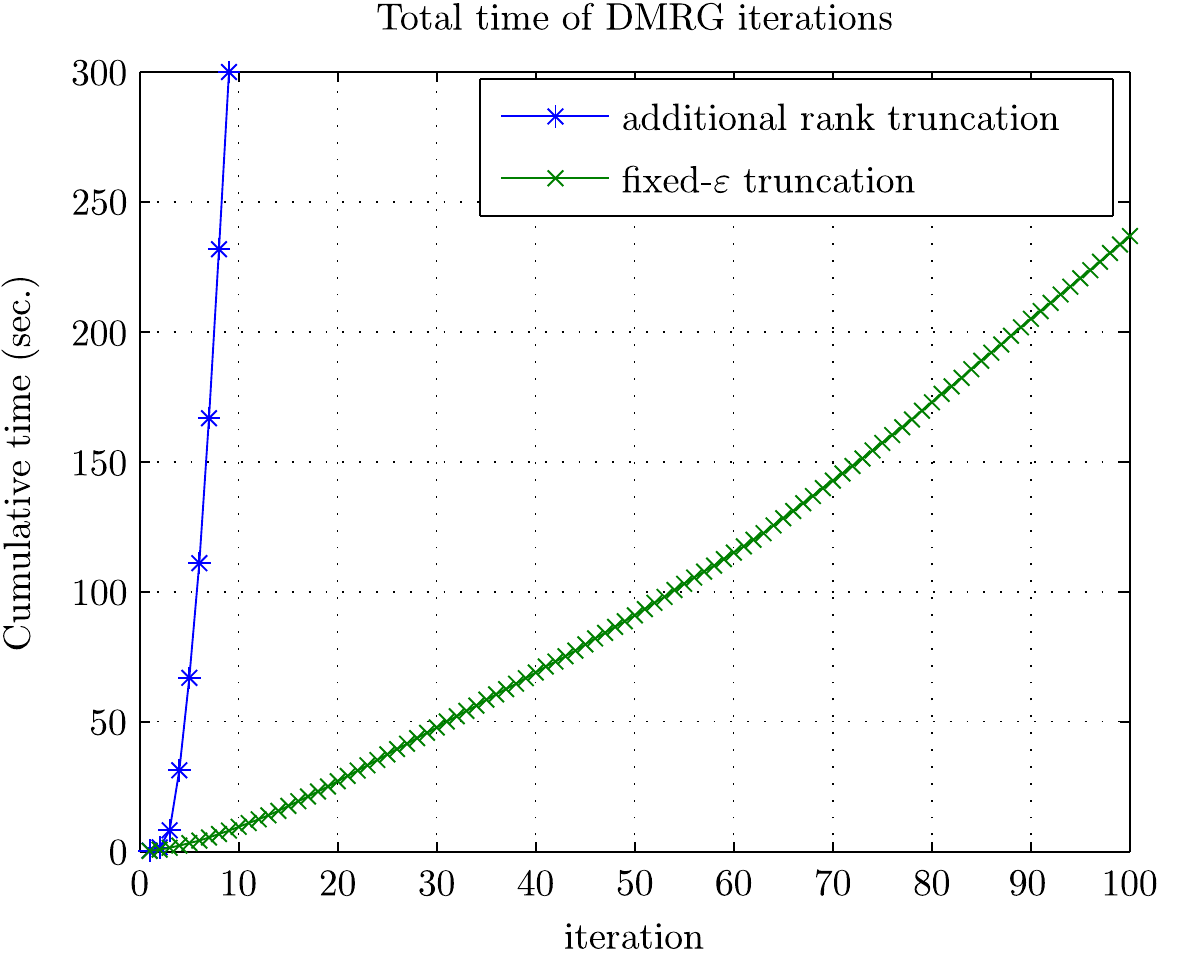}
 \caption{Cumulative times, sPDE, DMRG solvers}
 \label{fig:cumtime_dmrg}
\end{minipage}
\end{figure}

\subsection{1D stochastic (parametric) PDE (Tables \ref{tab:1dspde_d20}-\ref{tab:1dspde_eps}, Figures \ref{fig:resid_dmrg}, \ref{fig:cumtime_dmrg})}
In this example we consider a 1D stochastic (multiparametric) equation from \cite{khos-pde-2010}:
\begin{equation}
-\dfrac{\partial}{\partial x} a(x,{\bf{y}}) \dfrac{\partial u(x,{\bf{y}})}{\partial x} = f(x) = 1~~\mbox{in}~\Omega \times Y = [-1,1] \times [-1,1]^d,
\label{eqn:1dspde_problem}
\end{equation}
Dirichlet boundary conditions on $\partial\Omega$, and the coefficient is given as a Karhunen-Loeve expansion:
$$
a(x,{\bf{y}}) = a_0(x) + \sum\limits_{j=1}^d \sqrt{\lambda_j} a_j(x) y_j, ~\mbox{with}
$$
$$
a_0(x) = 1, \quad \sqrt{\lambda_j} = \dfrac{1}{2(j+1)^2}, \quad a_j(x)=\sin(\pi j x).
$$
The problem is then $d+1$-dimensional, and is not tractable in the full format.
It is again discretized using the FD scheme with the collocation method in the parameters on uniformly
distributed points. We use the preconditioner \cite{dkot-P2-2011}
$$
P_2 = \Delta^{-1} \Gamma(1/a) \Delta^{-1},
$$
where $\Gamma(a)$ is a stiffness matrix of the discretized elliptic operator \eqref{eqn:1dspde_problem} with the coefficient $a$. The parametric inversed Laplacian reads just $\Delta^{-1}_x \otimes I_{y_1} \otimes \cdots \otimes I_{y_d}$. Moreover,
we used the \emph{QTT} format in this example, with the explicit analytic QTT representation of the 1D
$\Delta^{-1}_x$ from \cite{khkaz-lap-2010}. To compute the reciprocal coefficient, we used the TT-structured Newton iterations. In the following, unless specially noted (table \ref{tab:1dspde_eps}), we fix the tensor rounding accuracy to $\eps = 10^{-5}$.

This example is essentially high-dimensional, with large ranks of the solution, and what is more important, of the coefficients. Hence we have to use the DMRG compression routines. The increased-rank truncation strategy allows to keep the accuracy, correspondingly increasing the time. But without it, one might get no relevant solution at all. We will demonstrate it in a comparison with the DMRG-solve algorithm.

We show in Tables \ref{tab:1dspde_d20}, \ref{tab:1dspde_d40}, \ref{tab:1dspde_d80} the number of iterations (it.), solution time (T, sec.), stabilized preconditioned residual (resid.) and the maximal TT rank versus the spacial $n_x$ and parametric $n_y$ grid sizes and the number of parameters $d$.

\begin{table}[!h]
\caption{$d=20$, $\eps=10^{-5}$}
\centering
\label{tab:1dspde_d20}
\begin{tabular}{|c|c||c|c|c|c|}\hline
 $n_x$ & $n_y$ & it. & T (sec.) & resid. & rank \\ \hline
 128   & 64   & 3 & 129.2 & 3.19e-6 & 28 \\ \hline
 256   & 64   & 3 & 124.1 & 2.93e-6 & 28 \\ \hline
 128   & 128  & 3  & 133.8 & 4.64e-6 & 27 \\ \hline
 128   & 256  & 3  & 148.3 & 4.68e-6 & 28 \\ \hline
\end{tabular}
\end{table}

\begin{table}[!h]
\caption{$d=40$, $\eps=10^{-5}$}
\centering
\label{tab:1dspde_d40}
\begin{tabular}{|c|c||c|c|c|c|}\hline
 $n_x$ & $n_y$ & it. & T (sec.) & resid. & rank \\ \hline
 128   & 64   & 3  & 413.7 & 2.13e-5 & 33 \\ \hline
 256   & 64   & 3  & 409.8 & 1.93e-5 & 33 \\ \hline
 128   & 128  & 3 & 334.7 & 1.51e-5 & 36 \\ \hline
 128   & 256  & 3 & 456.3 & 1.90e-5 & 33 \\ \hline
\end{tabular}
\end{table}

\begin{table}[!h]
\caption{$d=80$, $\eps=10^{-5}$}
\centering
\label{tab:1dspde_d80}
\begin{tabular}{|c|c||c|c|c|c|}\hline
 $n_x$ & $n_y$ & it. & T (sec.) & resid. & rank \\ \hline
 128   & 64   & 3 & 1187 & 1.71e-5 & 37 \\ \hline
 256   & 64   & 3  & 1280& 1.70e-5 & 36 \\ \hline
 128   & 128  & 3  & 1122 & 1.82e-5 & 35 \\ \hline
 128   & 256  & 3  & 1336 & 2.03e-5 & 33 \\ \hline
\end{tabular}
\end{table}

Consider a dependence on the tensor rounding accuracy $\eps$ in the case $n_x=128$, $n_y=64$, $d=20$. As in the previous tables,  we show the number of iterations, solution time, stabilized residual and maximal TT rank, see \ref{tab:1dspde_eps}.
\begin{table}[h!]
\caption{Dependence on the rounding accuracy $\eps$. Problem sizes $n_x=128$, $n_y=64$, $d=20$.}
\centering
\label{tab:1dspde_eps}
\begin{tabular}{|c||c|c|c|c|} \hline
$\eps$	& $10^{-3}$	& $10^{-4}$	& $10^{-5}$	& $10^{-6}$ \\ \hline
it.	& 2		& 2		& 3		& 3	\\ \hline
T, sec.	& 7.31		& 16.98		& 129.2		& 384.41	\\ \hline
resid.	& 2.05e-4	& 7.90e-5	& 3.19e-6	& 7.31e-8	\\ \hline
rank	& 8		& 13		& 28		& 46	\\ \hline
\end{tabular}
\end{table}
With increasing accuracy, the computational time increases drastically, as it depends both on the number of iterations and TT ranks.

Now, consider the TT-DMRG-solver from \cite{DoOs-dmrg-solve-2011} applied to the same problem $P_2 \Gamma(a) u = P_2 f$ with
$n_x=128$, $n_y=64$ and $d=20$. Following the section \ref{sec:dmrg_round} we compare two variants of rank truncations in the superblock splitting: with fixed $\eps$ and additional rank increasing. The convergence histories are shown on Fig. \ref{fig:resid_dmrg}, and the cumulative times on Fig. \ref{fig:cumtime_dmrg}, respectively.

We see, that increasing of truncation ranks improves the convergence significantly, despite that the time of each iteration is larger than with the fixed-accuracy truncation.

Notice the difference in time between the GMRES and DMRG solver. To achieve the same residual $4 \cdot 10^{-6}$ GMRES spent 129 sec., whereas DMRG (with increased ranks only) took about 300 sec. This shows the advantage of rapidly converging GMRES, provided a good preconditioner is given. It is natural that the also DMRG-based approximate MatVecs (which are in fact, just the DMRG truncations, up to additional structure of TT blocks, provided by their construction as Matrix by Vector multiplications) are cheaper than the linear system solutions.

\section{Conclusion}
The adapted GMRES method in the TT format for a linear system solution was proposed and investigated.
For the method presented the error analysis and performance improvements obtained with the aid of the inexact Krylov methods theory. The numerical experiments show, that the method provides a linear with respect to the grid size complexity in the case of TT approximation, and even logarithmic complexity with the QTT format.
The method was compared with the direct ALS/DMRG-type minimization solver for the TT format.
These methods manifest comparable timings and accuracies, and the GMRES method might be recommended in cases, when a good preconditioner is known.

\ifx\undefined\BibPrefix\def\BibPrefix{}\else\fi

\end{document}